\documentclass[a4paper,11pt]{article}
\usepackage{hyperref}
\hypersetup{
    colorlinks=true,
    linktoc=all,
    linkcolor=red,     
    citecolor=blue,     
}
\usepackage{graphicx}
\usepackage{amssymb}
\usepackage{latexsym, bm}
\usepackage{multicol}
\usepackage{indentfirst}
\usepackage{amsfonts}
\usepackage{amsmath}
\usepackage{setspace}
\newtheorem{theorem}{Theorem}[section]
\newtheorem{lemma}{Lemma}[section]
\newtheorem{claim}{Claim}[section]

\newtheorem{corollary}{Corollary}[section]

\newtheorem{definition}{Definition}[section]
\newtheorem{problem}{Problem}[section]

\newcommand{\ignore}[1]{}

\begin{document}
\begin{spacing}{0.98}

\title{Fan-complete Ramsey numbers}
\date{}

\author{
Fan Chung\footnote{Department of Mathematics, University of California at San Diego, United States of America. Email: {\tt fan@ucsd.edu}. } \;\;  and \;\; Qizhong Lin\footnote{Center for Discrete Mathematics, Fuzhou University,
Fuzhou, 350108, P.~R.~China. Email: {\tt linqizhong@fzu.edu.cn}. Supported in part by National Key R\&D Program of China (Grant No. 2023YFA1010202), NSFC (No.\ 12171088)  and NSFFJ (No.\ 2022J02018).}
}

\maketitle

\begin{abstract}
For graphs $G$ and $H$, we consider Ramsey numbers $r(G,H)$ with tight lower bounds, namely,
$
r(G,H) \geq (\chi(G)-1)(|H|-1)+1,
$
where $\chi(G)$ denotes the chromatic number of $G$ and $|H|$ denotes the number of vertices in $H$.
We say $H$ is $G$-good if the equality  holds.

Let $G+H$ be the join graph obtained from graphs $G$ and $H$ by  adding all edges between the  disjoint vertex sets of $G$ and $H$.
Let $nH$ denote the union graph of $n$ disjoint copies of $H$. We show that $K_1+nH$ is $K_p$-good if $n$ is sufficiently large.
In particular, the fan-graph $F_n=K_1 + n K_2$ is $K_p$-good if $n\geq 27p^2$, improving previous tower-type lower bounds for $n$ due to Li and Rousseau (1996). 
Moreover, we give a stronger lower bound inequality for Ramsey number $r(G, K_1+F)$ for the case of $G=K_p(a_1, a_2, \dots, a_p)$, the complete $p$-partite graph with $a_1=1$ and $a_i \leq a_{i+1}$.
In particular, using a stability-supersaturation lemma by Fox, He and Wigderson (2021), we show that for any fixed graph $H$,
\begin{align*}
r(G,K_1+nH) = \left\{ \begin{array}{ll}
(p-1)(n |H|+a_2-1)+1 & \textrm{if $n|H|+a_2-1$ is even or $a_2-1$ is even,}\\
(p-1)(n |H|+a_2-2)+1 & \textrm{otherwise,}
\end{array}
 \right.
\end{align*}
where $G=K_p(1,a_2, \dots, a_p)$ with $a_i$'s satisfying some mild conditions and $n$ is sufficiently large.
The special case of  $H=K_1$ gives an answer to Burr's question (1981) about the discrepancy of $r(G, K_{1,n})$ from $G$-goodness  for sufficiently large $n$. All bounds of $n$ we obtain are not of tower-types.

\medskip

\end{abstract}

\section{Introduction}

For graphs $G$ and $H$, the Ramsey number $r(G,H)$ is the smallest positive integer $N$ such that any graph  on $N$ vertices contains $G$ as a subgraph, or its complement contains $H$ as a subgraph.
A classic result of Chv\'{a}tal \cite{c} states
\begin{align}
r(K_p, T_n)=(p-1)(n-1)+1,
\label{chvatal}
\end{align}
where  $K_p$ is the complete graph on $p$ vertices and  $T_n$ is a tree with $n$ vertices.
Let $H$ be a connected graph with $n$ vertices. Since the graph consists of $p-1$ disjoint copies of $K_{n-1}$ is $H$-free and its complement is $K_p$-free,  one can easily derive (see \cite{c-h})
\begin{align}\label{c-h}
r(K_p,H)\ge (p-1)(n-1)+1.
\end{align}

For a graph $G$, let  $\chi(G)$ be the chromatic number of $G$,
and $s(G)$ the {\em chromatic surplus} of $G$,  i.e., the minimum size of a color class over all proper vertex-colorings of $G$ with $\chi(G)$ colors. e.g., $s(K_p)=1$.
Burr \cite{bur}  improved the lower bound  in (\ref{c-h}) by showing
\begin{align}\label{burr}
r(G,H)\ge(\chi(G)-1)(|H|-1)+s(G),
\end{align}
where  $H$ is a connected graph with $|H|\ge s(G)$ vertices.
A  graph $H$ is said to be {\em $G$-good} if the equality in (\ref{burr}) holds. For example,  all trees are $K_p$-good from the result of Chv\'{a}tal \cite{c}.

Burr and Erd\H{o}s \cite{be} initiated the study of Ramsey goodness problems that have since attracted the attention of many researchers.   For various generalizations of the goodness problems, the reader  is referred to the survey \cite{cfs-15} by Conlon, Fox and Sudakov.

For graphs $G$ and $H$, let $G+H$ be the join graph obtained from two graphs $G$ and $H$ by connecting the disjoint vertices of $G$ and $H$ completely.
Let $nH$ denote the union graph of $n$ disjoint copies of $H$. 
The fan graph $F_n$, also called the friendship graph, is the graph consisting of $n$ triangles all sharing a vertex, i.e., $F_n=K_1+nK_2$.

As an early application of the Erd\H{o}s-Simonovits stability lemma \cite{e3,er66,sim}, Li and Rousseau \cite{lr-jgt} showed that for any fixed graphs $G$ and $H$ and sufficiently large $n$,
\begin{align}\label{lr}
r(K_2+G,K_1+nH)=(\chi(G)+1)n|H|+1.
\end{align}
 This implies that $K_1+nH$ is $(K_2+G)$-good for sufficiently large $n$. In particular, $F_n$ is $K_p$-good for $p\ge3$ and sufficiently large $n$.
However,  the original stability results utilize a modified form of  progressive induction and, therefore, the lower bound for $n$ in (\ref{lr}) is quite large  as a form of tower type.

As a special case, the fan-complete Ramsey number has attracted much of attention.
In particular, it is known that $F_n$ is $K_p$-good for  $p=3,4,5,6$ and $n\ge p$, see \cite{lr-jgt,sbb,cz,kos}. In \cite{hht}, the authors claimed that $F_n$ is $K_p$-good if $n>cp^2$ for some constant $c>0$, but the paper contains a critical error
in \cite[Lemma 2.3]{hht} (lines 4--6, page 66,
 while using induction without enough  vertices in the neighborhood of a vertex)\footnote{The authors \cite{hht-pc} corrected their proof by showing that $F_n$ is $K_p$-good if $n>cp^3$ for some constant $c>0$.}.

 Let $B_{k,n}$ be the book graph on $n$ vertices which consists of $n-k$ copies of $K_{k+1}$ all sharing a common $K_k$, i.e. $B_{k,n}= K_k+(n-k)K_1$. 
Using the regularity lemma \cite{sz}, Nikiforov and Rousseau \cite{nr} showed that the book graph $B_{k,n}$ is $K_{p}$-good for fixed positive integers $k,p$ and sufficiently large $n$. Furthermore, extending the method used in \cite{nr}, Nikiforov and Rousseau  \cite{nr09} obtained a number of general goodness results.
However, all the bounds on $n$ of these results are of tower-types since the proofs rely on the regularity lemma.
Recently, using a stability-supersaturation lemma instead  of the regularity lemma, Fox, He and Wigderson \cite{fox} proved that
$B_{k,n}$ is $K_p$-good if $n\ge2^{k^{10p}}$.

\medskip
In this paper, we  first prove the following theorem whose proof can be found in Section \ref{pf-1}.
\begin{theorem}\label{main}
For any graph $H$ and $p\ge1$, $ r(K_p,K_1+nH)=(p-1)n|H|+1$ for $n\ge c{p\ell}/{|H|}$ where $\ell=r(K_p,H)$ and $c=(3+3\sqrt{2})^2 \approx 52.456$.
Namely,  $K_1+nH$ is $K_p$-good if $n\ge c{p\ell}/{|H|}$.
\end{theorem}

The following corollary is an immediate consequence of Theorem \ref{main}.
\begin{corollary}\label{fan}
If $n\ge Cp^2$,  then $F_n$ is $K_p$-good where $C= (3+3\sqrt{2})^2/2 \approx 26.228$.
\end{corollary}

The second part of this paper concerns improvements or generalizations of the lower bounds of Chv\'atal and Burr in (\ref{c-h}) and (\ref{burr}). 
We note that it is known  \cite{aks,kim} that  $r(K_3,K_n)=\Theta(n^2/\log n)$  while the lower bounds using (\ref{c-h}) and (\ref{burr}) only give $r(K_3,K_n)\ge 2(n-1)+1$.
A natural question is to
generalize the lower bound inequalities in (\ref{c-h}) and (\ref{burr})   and find families of graphs achieving or close to such lower bounds.  Nevertheless, there are some obvious obstacles.
For example,  the star $K_{1,n}$ is not $C_4$-good 
for all sufficiently large $n$ since
\[  n + \lceil \sqrt{n} \rceil +1 \geq r(C_4,K_{1,n})\geq  n+\lfloor n^{1/2}-6n^{11/40}\rfloor, \]
where the upper bound can be derived from the Tur\'an number of $C_4$ and the lower bound can be found in \cite{befrs} using probabilistic arguments.  F\"{u}redi \cite{furedi} proved that
 $r(C_4,K_{1,n})= n+\lceil n^{1/2} \rceil$, for infinitely many $n$.

For graphs $G$ and $H$, let us define the discrepancy of $r(G,H)$ from $G$-goodness as follows.

\begin{definition} Let $G$ be a graph  with chromatic surplus $s(G)$, and let $H$ be a connected graph  on at least $s(G)$ vertices.
Define $d(G,H)$ as the discrepancy of $r(G,H)$ from $G$-goodness, i.e., $d(G,H)=r(G,H)-(\chi(G)-1)(|H|-1)-s(G)$.
\end{definition}

From the above definition, $d(G,H)=0$ if and only if  $H$ is $G$-good. 
 Let  $K_p(a_1, \dots, a_p)$ denote  the complete $p$-partite graph  with parts of size $a_1,\dots,a_{p}$.
We know from Chv\'{a}tal and Harary \cite{c-h-1} that $r(K_{1,a_2},K_{1,n})=a_2+n-1$ if both $a_2$ and $n$ are even, and $r(K_{1,a_2},K_{1,n})=a_2+n$ otherwise. So $d(K_{1,a_2},K_{1,n})= a_2-2$ if both $a_2$ and $n$ are even, and $d(K_{1,a_2},K_{1,n})= a_2-1$  otherwise as observed by Burr  \cite{bur}. The discrepancy of $r(K_{1,a_2},K_{1,n})$ from $K_{1,a_2}$-goodness grows as $a_2$ grows. 
For $p \geq 3$, Burr \cite{bur} asked the question of determining when the {\it discrepancy} of $r(G, K_{1,n})$ from $G$-goodness grows in $n$ for $G=K_p(1,a_2, \dots, a_p)$.


To address the question of the discrepancy of $K_p(1,a_2, \dots, a_p)$ and $K_{1,n}$, we will first derive the following lower bound that improves the inequality of Burr in  (\ref{burr}) in some cases.

\begin{theorem}\label{low-bou}
Let $G=K_p(a_1, \dots, a_p)$ where $1= a_1\le a_2\leq \dots \leq a_{p}$. For any graph $F$ of order $n\ge2a_2$,
\begin{align*}
r(G,K_1+F) \ge \left\{ \begin{array}{ll}
(p-1)(n +a_2-1)+1 & \textrm{if $n+a_2-1$  is even or $a_2-1$ is even,}\\
(p-1)(n +a_2-2)+1 & \textrm{otherwise.}
\end{array} \right.
\end{align*}
\end{theorem}

Furthermore, we will show that the above lower bound is sharp in some general setting:
\begin{theorem}\label{main-2}
For any fixed graph $H$, integers $p\ge2$ and $b\ge1$, there exists $\delta>0$ such that the following holds for all  $n\ge pb^2/\delta$.
Let $a_1,\dots,a_p$ be positive integers  with $1= a_1\le a_2\leq \dots \leq a_{p-1}\leq b$ and $a_{p}\leq\delta n$, and let $G=K_p(1,a_2, \dots, a_p)$.
Then
\begin{align*}
r(G,K_1+nH) = \left\{ \begin{array}{ll}
(p-1)(n |H|+a_2-1)+1 & \textrm{if $n|H|+a_2-1$  is even or $a_2-1$ is even,}\\
(p-1)(n |H|+a_2-2)+1 & \textrm{otherwise.}
\end{array} \right.
\end{align*}

Moreover, we may take $\delta$ with $0<\delta<\min\big\{\frac{1}{400a^{|H|+2}p^4},(100a^pp^{14p})^{-A}\big\}$, where $a=\sum_{i=1}^{p-1}a_i$ and $A=\prod_{i=1}^{p-1}a_i$.
\end{theorem}

We remark that  the special case of  $H=K_1$ in Theorem \ref{main-2} gives an answer to Burr's question about the discrepancy
of $r(G, K_{1,n}) $ from $G$-goodness for sufficiently large $n$. In particular, the discrepancy of $r(G,K_{1,n})$ from $G$-goodness grows as $p$ and $a_2$ grow.

\medskip
The following corollary improves  (\ref{lr}) since  for any fixed graph $G$ with chromatic number $p$,  $K_2+G$ is a subgraph of  $K_{p+2}(1,1,a_3,\dots,a_{p})$ for some  $a_3,\dots,a_{p}$. Furthermore,  the lower bound on $n$ we obtain is not of tower-type.
\begin{corollary}\label{cor-lr}
For any fixed graph $H$, integers $p\ge2$ and $b\ge1$, there exists $\delta>0$ such that the following holds for all  $n\ge pb^2/\delta$. Let $a_1,\dots,a_p$ be positive integers  with $a_1\le a_2\leq \dots \leq a_{p-1}\leq b$ and $a_{p}\leq\delta n$.  If $a_1=a_2=1$, then $K_1+nH$ is $K_p(a_1, \dots, a_p)$-good.
\end{corollary}

In this paper, we use the following notation:  For a graph $G=(V,E)$ with vertex set $V$ and edge set $E$, we use $e(G)$ to denote the number of edges  $|E|$ in $G$. For $X \subseteq V$, let $G[X]$ denote the subgraph of $G$ induced by $X$, and let $e(X)$ denote the number of edges in $G[X]$.
For two disjoint subsets $X,Y\subseteq V$, we use $e(X,Y)$ to denote the number of edges between $X$ and $Y$.
For a vertex $v\in V$ and $X \subseteq V$, we denote by $N_X(v)$ the neighborhood of $v$ in $X$, and let $d_X(v)=|N_X(v)|$. The neighborhood of a vertex $v$ in $G$ is denoted by $N_G(v)$, i.e. $N_G(v)=N_V(v)$ and the degree of $v$ in $G$ is $d_G(v)=|N_G(v)|$.
$A \sqcup B$ denotes the disjoint union of $A$ and $B$. 
Let $[p]=\{1,2,\dots,p\}$. For undefined terminology, the reader is referred to \cite{bm}.

\section{Proof of Theorem \ref{main}}\label{pf-1}

Let $K_p(s)$ denote the complete $p$-partite graph $K_p(\underbrace{s,\dots,s}_{p})$.
In order to prove Theorem \ref{main}, we need the following lemma, which will also be applied in Theorem \ref{main-2}.
\begin{lemma}\label{turan}
For $ p\ge2$, let $\Gamma$ be a subgraph of $K_p(s)$. If $\Gamma$ has $z< s^2$ non-edges, then $\Gamma$ contains at least
$$
s^{p-2}(s^2-z)
$$
distinct copies of  $K_p$.
\end{lemma}
{\bf Proof.} Note that there are $s^p$ distinct copies of $K_p$ in $K_p(s)$ and each non-edge of $\Gamma$ destroys at most $s^{p-2}$ distinct copies of $K_p$. Therefore, if there are $z$ non-edges of $\Gamma$, then there are still at least $s^p-z\cdot s^{p-2}$ distinct copies of $K_p$ remaining.
 The proof is finished.
  \hfill  $\Box$

\medskip
We also have the following simple lemma.
\begin{lemma}\label{deg} For any graph $H$ and integers $p,n\ge1$,
$r(K_p,nH)\le h(n-1) +r(K_p,H).$
\end{lemma}
{\bf Proof.} The base case where $n=1$ is clear. Suppose the assertion holds for $n-1$, then any graph $G$ of order $h(n-1) +r(K_p,H)$ contains a $K_p$ or its complement $\overline{G}$ contains $(n-1)H$. Suppose there is no $K_p$, otherwise we are done. Thus $\overline{G}$ contains $(n-1)H$. Deleting the $(n-1)h$ vertices of $(n-1)H$ from $G$, there are $r(K_p,H)$ vertices remaining. Let $G'$ be the subgraph induced by these remaining vertices. $\overline{G'}$ must contain a copy of $H$ since $G'$ contains no $K_p$ from the assumption. However then together with the previous $(n-1)H$ will yield $nH$ in $\overline{G}$. \hfill$\Box$

\medskip
\noindent
{\bf Proof of Theorem \ref{main}.} The lower bound follows from (\ref{burr}), so we will focus on the upper bound in the following.
Let $N=(p-1)hn+1$, where $h=|H|$ is the order of the graph $H$  and  $n\ge (3+3\sqrt{2})^2{p\ell}/{h}$ with $\ell=r(K_p,H)$.
The assertion is clear for $p=1,2,$ so we may assume $p\ge3$.
Suppose to the contrary that  there exists a graph $\Gamma$ on $N$ vertices such that $\Gamma$ is $K_p$-free and its complement $\overline{\Gamma}$ contains no copy of $K_1+nH$. We shall show that this leads to a contradiction.

Let $V$ denote the vertex set of $\Gamma$. For any vertex $v\in V$, by Lemma \ref{deg}, we have 
\begin{align}\label{deg-upp}
d_{\overline{\Gamma}}(v)<r(K_p,nH)\le h(n-1) +\ell.
\end{align}
It follows that
\begin{align}\label{deg-low}
d_\Gamma(v)= N-1-d_{\overline{\Gamma}}(v)\ge N-h(n-1) -\ell\ge(p-2)hn+h-\ell.
\end{align}

\begin{claim}\label{inte-up}
There exists a partition $\sqcup_{i=1}^{p-1}V_i$ of $V$ such that the total number of internal edges is at most $\frac{N}{2}(\ell-h)$, i.e., $\sum_{i=1}^{p-1}e(V_i)\le \frac{N}{2}(\ell-h)$.
\end{claim}
{\bf Proof.}
We apply the degree majorization algorithm used by Erd\H{o}s \cite{e70} and  F\"{u}redi \cite{fur}.
Let $V_{0}^+=V$. For $i\ge1$, pick a vertex $v_i\in V_{i-1}^+$ such that $v_i$ has the maximum degree in $\Gamma[V_{i-1}^+]$, and let $V_i=V_{i-1}^+\setminus N_\Gamma(v_i)$ and $V_{i}^+=V_{i-1}^+\cap N_\Gamma(v_i)$.
The procedure stops  when there are no other vertices remaining. Let $r$ be the largest $i$ such that $V_i$ is defined. Clearly, $V_1,\dots,V_{r-1},V_{r}=V_{r-1}^+$ form a partition of $V$.
Note that $r\le p-1$ since $\{v_1,\dots,v_r\}$ induces a complete graph. Combining $|V_i|\le h(n-1) +\ell$ from (\ref{deg-upp}) and  $(p-2)(h(n-1) +\ell)<N$,
we have $r\ge p-1$.
Therefore, we conclude $r=p-1$.

Note that for $i\in[p-1]$ and $x\in V_i$, $d_{V_{i-1}^+}(x)\le d_{V_{i-1}^+}(v_i)=|V_{i}^+|$ from the choice of vertex $v_i$. Thus we have $2e(V_i)+e(V_i,V_{i}^+)=\sum_{x\in V_i}d_{V_{i-1}^+}(x)\le|V_i||V_{i}^+|.$
By adding up both sides of the inequality, we have
\begin{align*}
e(\Gamma)+\sum_{i=1}^{p-1}e(V_i)\le \sum_{i=1}^{p-1}|V_i||V_{i}^+| \leq e(T_{N,p-1}),
\end{align*}
 where $T_{N,p-1}$ denotes the Tur\'an graph on $N$ vertices containing no $K_p$ with the maximum number of edges.
Combining with the fact that  $e(\Gamma)\ge\frac{N}{2}(N-h(n-1) -\ell)$ from (\ref{deg-low}), we conclude that the total number of internal edges
satisfies
\begin{align*}
\sum_{i=1}^{p-1}e(V_i)\le \left(1-\frac{1}{p-1}\right)\frac{N^2}{2}-\frac{N}{2}(N-h(n-1) -\ell)
\le\frac{N}{2}(\ell-h),
\end{align*}
as claimed.\hfill$\Box$

\medskip 

Now, let us take a partition $\sqcup_{i=1}^{p-1}V_i$ such that it attains the minimal number $\sum_{i=1}^{p-1}e(V_i)$ of the internal edges. Thus, we must have that
\begin{align}\label{mo-dg-oth}
\text{for each vertex $v\in V_i$ and $j \not = i$, $d_{V_i}(v)\leq d_{V_j}(v)$,}
\end{align}
since otherwise there exists a vertex $v\in V_i$ with $d_{V_i}(v)> d_{V_j}(v)$ for some $j\neq i$, and we can then  put $v$ into $V_j$ to get a smaller total number of internal edges since the number of the internal edges of the `new partition' will decease by $d_{V_i}(v)-d_{V_j}(v)$,
 which is a contradiction.

We set  $m=\frac{N}{2}(\ell-h)$. 
From  (\ref{deg-low}) and Claim \ref{inte-up}, we  have
\begin{align}\label{edge-bt}
\sum_{1\le i<j\le p-1} e(V_i,V_j)
&= e(\Gamma) -\sum_{i=1}^{p-1} e(V_i)  \nonumber \\
&\ge\frac{N}{2}(N-h(n-1) -\ell)-m
\nonumber \\
&\ge\left(1-\frac{1}{p-1}\right)\frac{N^2}{2}-2m.
\end{align}
Using $\sum_{i=1}^{p-1} (|V_i|-\frac{N}{p-1})^{2}=\sum_{i=1}^{p-1} |V_i|^{2}-\frac{N^{2}}{p-1},$
 we have, for each $i\in[p-1]$,
\begin{align}\label{Vi-b}
\left||V_i|-\frac{N}{p-1}\right|\le2\sqrt{m},
\end{align}
since otherwise $\sum_{i=1}^{p-1} |V_i|^{2}>\frac{N^{2}}{p-1}+4m$ and so
\begin{align*}
\sum_{1\le i<j\le p-1} e(V_i,V_j)\le\sum_{1\le i<j\le p-1} |V_i||V_j|=\frac{1}{2}\left(N^{2}-\sum_{i=1}^{p-1} |V_i|^{2}\right)
< \frac{1}{2}\left(N^{2}-\frac{N^{2}}{p-1}\right)-2m,
\end{align*}
which contradicts (\ref{edge-bt}).

For distinct $i$ and $j$, $1\leq i,j\leq p-1$,
we define $z_{i,j} = |V_i||V_j|-e(V_i,V_j)$, which is the number of non-edges between $V_i$ and $V_j$.
Let $z=\sum_{i<j} z_{i,j}$. Then we must have
\begin{align}\label{Vi-Vj}
z \leq 2m
\end{align}
since
$\sum_{1\le i<j\le p-1} |V_i||V_j|-z=\sum_{1\le i<j\le p-1} e(V_i,V_j)\geq (1-\frac{1}{p-1})\frac{N^{2}}{2}-2m$
from (\ref{edge-bt}) and the fact that $\sum_{1\le i<j\le p-1} |V_i||V_j|$ is at most the Tur\'{a}n number of $K_p$-free graph on $N$ vertices.

In the following, we will show that for each vertex $v\in V_i$, $d_{V_i}(v)\leq \sqrt{2m}$.
Suppose to the contrary that without loss of generality there exists some vertex $v\in V_1$  having $s>\sqrt{2m}$ neighbors in its own part $V_1$. It follows from (\ref{mo-dg-oth}) that $v$ also has at least $s$ neighbors in each of the other parts. Let $U_i=N_\Gamma(v)\cap V_i$ denote the neighborhood of $v$ in $V_i$ for $i=1,2,\dots,p-1$. Clearly, $|U_i|\ge s>\sqrt{2m}$. It follows by Lemma \ref{turan} that there are at least
$$
s^{p-3}(s^2-z)\ge s^{p-3}(s^2-2m)
$$
copies of $K_{p-1}$ in the neighborhood of $v$. Therefore, $\Gamma$ definitely contains a copy of $K_p$, which leads to a contradiction.

Let $t=(2+\sqrt{2})\sqrt{m}+\ell-h$. We claim that
\begin{align}\label{Vi-out-d}
\text{for each vertex $v\in V_i$ and $j\neq i$,  $d_{V_j}(v)\geq|V_j|-t$.}
\end{align}
On contrary, suppose that some vertex $v\in V_i$ has at least $t$ non-neighbors in $V_j$. From the above, $v$ has at least $|V_i|-\sqrt{2m}-1$ non-neighbors in $V_i$. Note that $|V_i|\ge\frac{N}{p-1}-2\sqrt{m}$ from (\ref{Vi-b}). In total, the number of non-neighbors of $v$ is at least
\[
t+|V_i|-\sqrt{2m}-1\ge\frac{N}{p-1}+\ell-h>h(n-1) +\ell.
\]
This contradicts (\ref{deg-upp}) that $d_{\overline{\Gamma}}(v)\le h(n-1) +\ell-1$.

Suppose that there exists an edge $uv\in V_1$. Let $W_i$ denote the common neighborhood of $u$ and $v$ in $V_i$ for $2\le i\le p-1$.
Then, from (\ref{Vi-out-d}) and $|V_i|> hn-2\sqrt{m}$, we have that for each $2\le i\le p-1$,
\[
|W_i|\ge |V_i|-2t\ge hn-(6+2\sqrt{2})\sqrt{m}-2\ell+2h>\sqrt{2m},
\]
where the last inequality follows from $m=\frac{N}{2}(\ell-h)<\frac{1}{2}(p-1)\ell hn$, and $n\ge(3+3\sqrt{2})^2{p\ell}/{h}$ by the assumption.
Therefore, we can apply Lemma \ref{turan} again to get a copy $K_{p-2}$ in the common neighborhood of $u$ and $v$,
which  leads to a contradiction.
Consequently, $V_1$ forms an independent set. Similarly, $V_i$ forms an independent set for each $2\le i\le p-1$.

Now, on the average, there is some part $V_i$ of size at least $\left\lceil N/(p-1)\right\rceil=hn+1$ which forms an independent set from the above. Therefore, we can definitely get a copy of $K_1+nH$ in the complement of $\Gamma$. The final contradiction completes the proof of Theorem \ref{main}.
\hfill$\Box$

\section{Improving  lower bounds}\label{pf-31}

The lower bound for $r(G, K_1+F)$  in Theorem \ref{low-bou} is by construction.
We first have the following lemma.
\begin{lemma}\label{const} 
Suppose $n\ge 2a_2$. If $n+a_2-1$ is even or $a_2-1$ is even, then there exists an $(a_2-1)$-regular triangle-free graph of order $n+a_2-1$. If both $n+a_2-1$ and $a_2-1$ are odd, then there exists an $(a_2-1)$-regular triangle-free graph of order $n+a_2-2$.
\end{lemma}
{\bf Proof.}
We consider the following two cases:

\medskip\noindent
{\it Case 1:}  $n+a_2-1$ is even.

\medskip
Let $X$ and $Y$ be two sets of vertices of size $\lambda=(n+a_2-1)/2$, say,
 $X=\{x_1,\dots,x_\lambda\}$ and $Y=\{y_1,\dots,y_\lambda\}$. Let $\Lambda=\Lambda(X,Y)$ be a bipartite graph with two parts $X$ and $Y$, in which $x_k$ is adjacent to $y_\ell$ if and only if $\ell=k+i\pmod \lambda$ for $i=0,\dots,a_2-2$.
So $\Lambda$ is an $(a_2-1)$-regular $K_3$-free graph as desired since $n\ge a_2-1$.

\medskip\noindent
{\it Case 2:} $a_2-1$ is even.

\medskip 
For this case, we consider the following construction by Sidorenko \cite{sid} for solving a problem of Erd\H{o}s (see \cite{rr}).
Let $k=(a_2-1)/2$ and $\mu=n+a_2-1$, and let $\Lambda$ be the graph whose vertex set is $\mathbb{Z}_\mu$, where any two
vertices $i, j\in \mathbb{Z}_\mu$ are connected by an edge if and only if
\[
(i -j)\in \{\pm k, \pm(k + 1),\dots, \pm(2k - 1)\}.
\]
Since  $n\ge 2a_2> 2a_2-4$ and so $\mu\ge 6k-2$,
it follows that $\Lambda$ is an $(a_2-1)$-regular $K_3$-free graph of order $\mu$ as desired.
If $n< 2a_2-4$, then such a graph $\Lambda$ constructed as above may contain a triangle (e.g., when $a_2=3$ and $n=1$).

\medskip\noindent
{\it Case 3:}
If both $n+a_2-1$ and $a_2-1$ are odd, then  $n+a_2-2$ is even. Then similar to Case 1 we can construct an $(a_2-1)$-regular triangle-free graph on $n+a_2-2$ vertices.
 \hfill$\Box$

\medskip\noindent
{\bf Proof of Theorem \ref{low-bou}.} First, suppose that $n+a_2-1$ is even or $a_2-1$ is even. We will show $r(G,K_1+F)>(p-1)(n+a_2-1)$. From Lemma \ref{const}, there exists an $(a_2-1)$-regular triangle-free graph $\Lambda$ of order $n+a_2-1$. Let $\Lambda_i$, $i=1,2,\dots,p-1$, be disjoint copies of $\Lambda$ with vertex sets $V_i$. Let $\Gamma$ be the graph obtained from $\sqcup_{i=1}^{p-1} \Lambda_i$ by adding all edges between $V_i$ and $V_j$ for $1\le i<j\le p-1$. Then the complement of $\Gamma$ contains no $K_1+F$ since $\overline{\Gamma}$ is $(n-1)$-regular. In the following, we will prove that $\Gamma$ contains no $K_1+K_{p-1}(a_2)$ by induction on $p\ge2$.

Let $V=\sqcup_{i=1}^{p-1} V_i$.
It is true for $p=2$  since $\Lambda$ is $(a_2-1)$-regular. So we may assume that $p\ge3$ and the assertion holds for smaller $p$.
 Suppose to the contrary that $\Gamma$ contains a subgraph $K_1 + K_{p-1}(a_2)$.  The vertex set of $K_1+K_{p-1}(a_2)$ is denoted by $\{u_0\} \sqcup(\sqcup_{i=1}^{p-1} U_i$).
We can relabel the $V_i$'s so that  $u_0$ is in $V_{p-1}$. Since $\Lambda$ is $K_3$-free, $V_{p-1}$ can only contain vertices in at most
one of the $U_i$'s.  Furthermore, since $\sqcup_{i=1}^{p-2} \Lambda_i$ can not contain $K_1 + K_{p-2}(a_2)$ by the inductive hypothesis,  $V_{p-1}$ must contain some vertices in $\sqcup_{i=1}^{p-1} U_i$. Let $U_{p-1}$ denote the set with $U_{p-1} \cap V_{p-1} \neq \emptyset$.
Since $\Lambda$ is $(a_2-1)$-regular, then there exists a vertex in $U_{p-1}$ not in $V_{p-1}$.
Moreover, we have  $\sqcup_{i=1}^{p-2} U_i \subseteq V \setminus V_{p-1}$.
This guarantees a copy of $K_1 + K_{p-2}(a_2)$ in  $V \setminus V_{p-1}$, which contradicts the inductive hypothesis.
Therefore,  $\Gamma$ contains no $K_1+K_{p-1}(a_2)$ as claimed.

If both $n+a_2-1$ and $a_2-1$ are odd, then from Lemma \ref{const}, there exists an $(a_2-1)$-regular triangle-free graph of order $n+a_2-2$. Therefore, by a similar argument as above, we can obtain that $r(G,K_1+F)>(p-1)(n+a_2-2)$.
 This completes the proof of Theorem \ref{low-bou}.
\hfill$\Box$

\medskip
Let $G=K_p(a_1,\dots,a_{p})$, where $a_1\le a_2\leq \dots \leq a_{p}$, and let $G_1=K_{p-1}(a_1,\dots,a_{p-1})$.
A result of Burr \cite[Theorem 5]{bur} states that for any connected graph $H$ on $n$ vertices,
\begin{align}\label{b-it}
r(G,H)\ge r(G_1,H)+n-1.
\end{align}

The following can be viewed as a slight improvement of  (\ref{b-it}) for the case of  $H=K_1+F$.
\begin{corollary}\label{imp-lw}
Let $a_1,\dots,a_p$ be integers  with $1= a_1\le a_2\leq \dots \leq a_{p}$.
Let $G=K_p(a_1,\dots,a_{p})$ and $G_1=K_{p-1}(a_1,\dots,a_{p-1})$. If $F$ is a graph of order $n$ for some integer $n\ge2a_2$, then
\begin{align*}
r(G,K_1+F) \ge \left\{ \begin{array}{ll}
r(G_1,K_1+F)+n+a_2-1 & \textrm{if $n+a_2-1$ is even,}\\
r(G_1,K_1+F)+n+a_2-2 & \textrm{otherwise.}
\end{array} \right.
\end{align*}

\end{corollary}
{\bf Proof.}  We first consider the case when $n+a_2-1$ is even.
Let $\Gamma_1$ be a graph on $r(G_1,K_1+F)-1$ vertices such that $\Gamma_1$ does not contain  $G_1$ and the complement of $G_1$ does not contain  $(K_1+F)$.
Let $\Gamma_2$ be an $(a_2-1)$-regular bipartite graph of order $n+a_2-1$ as constructed in Lemma \ref{const}. It is not difficult to verify that $\Gamma_2$ is $K_{s,t}$-free for any $s,t$ with $s+t\ge a_2+1$ since $n\ge2a_2$.
Let $\Gamma$ be the join graph $\Gamma_1+\Gamma_2$. Clearly, $\overline{\Gamma}$ is $(K_1+F)$-free because its maximum degree is $n-1$.

We want to show that  $\Gamma$ contains no copy of $G=K_p(1,a_2,\dots,a_{p})$. Suppose that, on the  contrary,  $\Gamma$ contains $G$ as a subgraph. Let  the vertex set of $G$ be denoted as the disjoint union of $V_1, V_2, \dots, V_p$ where $|V_i|=a_i,
V_1= \{u\}$ and any edge of $G$ is between $V_i$ and $V_j$ for some $i \not = j$.

\medskip
{\it Case 1:}  $v$ is in $\Gamma_2$.

\medskip
 Since $\Gamma_2$ is $K_{1,a_2}$-free and $\Gamma_2$ is bipartite, at most one of the $V_i$'s contains vertices in $\Gamma_2$.
 Therefore, $\Gamma_1$ must contains a copy of $K_p(1,a_2,\dots,a_{p-1})$. This leads to a contradiction of the inductive assumption on $p$.

\medskip
 {\it Case 2:}  $v$ is in $\Gamma_1$.

\medskip
Since $\Gamma_2$ is bipartite, at most two  of the $V_i$'s  contain vertices in $\Gamma_2$.
Since $\Gamma_2$ is $K_{s,t}$-free for any $s,t$ with $s+t\ge a_2+1$, $\Gamma_1$ must contains a copy of $K_p(1,a_2-s',a_3-t'+1,\dots,a_{p})$ with $s'+t' \le a_2$. Note that $a_2-s'+a_3-t'+1\ge a_3$, which implies that $\Gamma_1$ contains a copy of $K_p(1,a_2,\dots,a_{p-1})$. This again leads to a contradiction.

\medskip
Therefore we conclude that $\Gamma$ contains no copy of $G=K_p(1,a_2,\dots,a_{p})$.

\medskip
For the case of $n+a_2-1$ odd, we use the same construction of $\Gamma_2$ with $n+a_2-2$ vertices satisfying the conditions
that $\Gamma_2$ is bipartite and $(a_2-1)$-regular. The proof is similar and will be omitted here.
\hfill$\Box$

\section{Proof of Theorem \ref{main-2}}\label{pf-2}

The following lemma, which is essentially due to Erd\H{o}s \cite{erdos64},
 states that if a graph $\Gamma$ on $N$ vertices  contains $\Omega(N^p)$ copies of $K_p$, then one can find a copy of $K_p(a_1,\dots,a_{p})$ in $\Gamma$ with one part  of size linear in $N$.
We here give a proof with specified bounds for various parameters, which will be used later in the proof of Theorem \ref{main-2}. The methods of the proof are similar to those in   \cite[Lemma 4.2]{fox}.
\begin{lemma}[Erd\H{o}s \cite{erdos64}]\label{Erdos}
For any $0<\delta<1$, and for any integers $p\ge2$, $b\geq1$, and $1\le a_1\le  \dots \leq a_{p-1}\leq b$, there exists some $\eta>0$ such that the following holds for all large $N$. If $\Gamma$ is a $K_p(a_1,\dots,a_{p})$-free graph on $N$ vertices with $a_p\leq \delta N$, then $\Gamma$ has at most $\eta N^p$ copies of $K_p$.

\medskip
Moreover, we may take $\eta=\delta^{1/(a_1a_2\cdots a_{p-1})}$, and $N\ge \max\{p^2,b^2/\delta\}$.
\end{lemma}
{\bf Proof.}
We define $\eta_s$, for $s=1,\dots, p-1$, inductively by choosing $\eta_p=\delta$, $\eta_{p-s}= (\eta_{p-s+1})^{1/a_{p-s}}$ and finally $\eta=\eta_1$. We will  show by induction on $s$ that a $K_{s+1}(a_{p-s},\dots,a_{p-1},a_p)$-free graph $\Gamma$ contains at most $\eta_{p-s} N^{s+1}$ copies of $K_{s+1}$.

For the base case of $s=1$, we want to show that a $K_2(a_{p-1},a_p)$-free graph $\Gamma$ contains at most $\eta_{p-1} N^2$ edges where $\eta_{p-1}= \delta^{1/a_{p-1}}$.  If $a_{p-1}=1$, then $e(\Gamma)\le a_pN/2<\eta_{p-1} N^2$. So the assertion obviously holds provided $N\ge a_p/\delta$. Thus we may assume $a_{p-1}\ge2$.  If $2e(\Gamma)/N\le a_{p-1}^2$, then we are done since $e(\Gamma)< a_{p-1}^2N\le b^2N\le \delta N^2\le \delta^{1/a_{p-1}} N^2$. Thus we assume $2e(\Gamma)/N\ge a_{p-1}^2$. We apply the double-counting method.
Since
 any $a_{p-1}$ vertices have at most $a_p\le\delta N$ common neighbors,  there are at most $\delta N{\binom N {a_{p-1}}}<\frac{\delta N^{a_{p-1}+1}}{a_{p-1}!}$ copies of $K_{1,a_{p-1}}$. Moreover, a vertex of degree $d$ contributes $\binom d {a_{p-1}}$ copies of $K_{1,a_{p-1}}$.
 Therefore, we have
 \[
\frac{\delta N^{a_{p-1}+1}}{a_{p-1}!}>\sum_{v\in V}{d_{\Gamma}(v)\choose a_{p-1}}\ge N{2e(\Gamma)/N\choose a_{p-1}}\ge \frac N{a_{p-1}!e}\left(\frac{2e(\Gamma)}{N}\right)^{a_{p-1}}
\]
by Jensen's inequality and the fact that
 ${t\choose p}\ge \frac{t^p}{p!e}$ for $t\ge p^2$ (since $(1-1/p)^{p-1}$ monotonically decreases to $1/e$).
Therefore, $e(\Gamma)<\delta^{1/a_{p-1}}N^2$, and  we may take $\eta_{p-1}= \delta^{1/a_{p-1}}=\eta_{p}^{1/a_{p-1}}$.

Now suppose $s\ge2$ and the assertion holds for any $s'< s$. We will show that any
$K_{s+1}(a_{p-s},\dots,a_{p-1},a_p)$-free graph $\Gamma$ contains at most $\eta_{p-s} N^{s+1}$ copies of $K_{s+1}$. Suppose to the contrary that $\Gamma$ contains at least $\eta_{p-s} N^{s+1}$ copies of $K_{s+1}$. For every $s$-set $S$ of $V(\Gamma)$, let ext($S$) be the set of vertices $v$ such that $S\cup\{v\}$ forms a $K_{s+1}$ in $\Gamma$. Note that the sum of $|{\rm ext}(S)|$ over all $s$-sets $S$ is exactly $s+1$ times the number of $K_{s+1}$ in $\Gamma$.  By the assumption, this sum is therefore more than $(s+1)\eta_{p-s} N^{s+1}$. Thus, the average value of $|{\rm ext}(S)|$ is greater than $(s+1)\eta_{p-s} N^{s+1}/{N\choose s}>(s+1)!\eta_{p-s} N$. Again by  Jensen's inequality, we have
\[
\sum_{S\subset V: |S|=s}{{\rm ext}(S)\choose a_{p-s}}\ge {N\choose s}{(s+1)!\eta_{p-s} N\choose a_{p-s}}.
\]
If $a_{p-s}=1$, then ${N\choose s}{(s+1)!\eta_{p-s} N\choose a_{p-s}}/{ N\choose a_{p-s}}={N\choose s}(s+1)!\eta_{p-s}\ge \eta_{p-s}N^s$ provided $N\ge s^2$. For $a_{p-s}\ge2$,  we choose $N\ge \max\{s^2,a_{p-s}^2/[(s+1)!\eta_{p-s}]\}$, and we have
\begin{align*}
{N\choose s}{(s+1)!\eta_{p-s} N\choose a_{p-s}}\bigg/{ N\choose a_{p-s}}\ge \frac{N^s}{s!e}\cdot\frac{[(s+1)!\eta_{p-s}N]^{a_{p-s}}}{a_{p-s}!e}\cdot\frac{a_{p-s}!}{N^{a_{p-s}}}
>\eta_{p-s}^{a_{p-s}}N^s.
\end{align*}
Consequently, we conclude that there is some $a_{p-s}$-set $T$ such that the common neighborhood of $T$ has more than $\eta_{p-s}^{a_{p-s}}N^s=\eta_{p-s+1}N^s$ copies of $K_s$. By the inductive assumption,  there is a copy of $K_{s}(a_{p-s+1},\dots,a_{p})$ in the common neighborhood of $T$. Together  with $T$ yields a copy of $K_{s+1}(a_{p-s},\dots,a_{p})$ in $\Gamma$, which leads to a contradiction. Therefore we have shown that a $K_{s+1}(a_{p-s},\dots,a_{p-1},a_p)$-free graph $\Gamma$ contains at most $\eta_{p-s} N^{s+1}$ copies of $K_{s+1}$ for  $N\ge \max\{s^2,a_{p-s}^2/[(s+1)!\eta_{p-s}]\}$ and
$$\eta:=\eta_1=\eta_2^{1/a_1}=\cdots=\eta_{p}^{1/(a_1a_2\cdots a_{p-1})}=\delta^{1/(a_1a_2\cdots a_{p-1})}.$$
This completes the proof of  Lemma \ref{Erdos}. \hfill$\Box$

\ignore{Based on the proof by Fox et al. \cite{fox}, the constant $\delta$  is approximately a double-exponential
in $\eta$, $p$ and $b$ (with base $\eta$), and  $N\ge  2pb/\eta$ will suffice.}
\medskip

We will apply the following {\em stability-supersaturation lemma} by Fox, He and Wigderson \cite[Theorem 3.1]{fox} (in a slightly different form) to obtain the desired structures for graphs forbidding some special classes of graphs.
Similar approaches are often referred to as  combinations of the stability theorem \cite{e3,er66,sim} and  the supersaturation result \cite{es}. This stability-supersaturation lemma implies that if a graph $\Gamma$ has slightly smaller minimum degree than the $K_p$-free Tur\'{a}n graph and has few copies of $K_p$, then it is close to the Tur\'{a}n graph.

\begin{lemma}[Fox, He and Wigderson \cite{fox}]\label{stability}
For every $\varepsilon>0$ and every integer $p\geq 2$, there exist $\eta,\gamma>0$ such that the following holds for all $N\geq 10$. Suppose $\Gamma$ is a graph on $N$ vertices with minimum degree at least $(1-\frac{1}{p-1}-\gamma)N$ and at most $\eta N^p$ copies of $K_p$. Then there is a partition $V(\Gamma)=\sqcup_{1\le i\le p-1} V_i$ such that the following hold:

\medskip
(i) $\sum_{1\le i\le p-1}e(V_i)\le\varepsilon {N\choose 2}$.

(ii) $\big||V_i|-\frac{N}{p-1}\big|\leq\sqrt{2\varepsilon}N$.

(iii) $e(V_i,V_j)\geq (1-p^2\varepsilon)|V_i||V_j|$.

(iv) For each $v\in V_i$, $d_{V_i}(v)\leq d_{V_j}(v)$.

\medskip
Moreover, we may take $\gamma=\min\{\frac{1}{2p^2},\frac{\varepsilon}{2}\}$ and $\eta=p^{-10p}\varepsilon$.
\end{lemma}

Note that   Conlon, Fox and Sudakov \cite[Corollary 3.4]{cfs}  obtained a stronger result by
 using  the minimum degree condition instead of the average degree condition and by using the  graph removal lemma (see e.g. Conlon and Fox \cite{cf}),  which, however, requires tower-type bounds in the parameters.

\medskip
Fox, He and Wigderson \cite{fox} established that if $n\ge2^{k^{10p}}$, then $B_{k,n}$ is $K_p$-good.
In order to give a better lower bound for $n$ of Theorem \ref{main-2},
we will use the following upper bound concerning the book graph $B_{k,n}$.
\begin{lemma}\label{fox}
Let $p,k,t\ge1$ be integers. Then $r(K_p,B_{k,t})\le (k+1)^pt$.
\end{lemma}
{\bf Proof.} The proof is by induction on $p\ge1$. The assertion is trivial for $p=1,2$, and so we may assume that $p\ge3$ and the assertion holds for smaller $p$.
Let $N_p=r(K_p,B_{k,t})-1$, and we consider a graph $\Gamma$ on $N_p$ vertices which contains no $B_{k,t}$ and its complement $\overline{\Gamma}$ is $K_p$-free. Let $V$ be the vertex set of $\Gamma$. By induction, we have $d_{\overline{\Gamma}} (v)\le N_{p-1}$ for any vertex $v\in V$. Thus, each vertex $v\in V$ has degree at least $N_p-N_{p-1}-1$ in $\Gamma$.

We first take an arbitrary vertex $v_1\in V$, and then we choose a neighbor, say $v_2$, of $v_1$ in $\Gamma$. Inductively, we can choose $k$ vertices $v_1,\dots, v_{k}$ which form a clique in $\Gamma$ and the number of the common neighbors of $v_1,\dots, v_{k}$ is at least $N_p-k(N_{p-1}+1)$. Since $\Gamma$ contains no $B_{k,t}$, we have
\[
N_p-k(N_{p-1}+1)< t-k.
\]
Therefore, it is not difficult to obtain that $N_p<(k+1)^pt$, completing the proof. \hfill$\Box$

\ignore{\begin{align*}
N_p&\le kN_{p-1}+ n-1
\\&\le k(kN_{p-2}+n-1)+ n-1
\\&\le k^{p-2}(kN_{1}+n-1)+k^{p-2}(n-1)+\cdots+k(n-1)+n-1
\\&< k^pn,
\end{align*}
completing the proof.
\hfill$\Box$}

\medskip
We remark  that if $p\ge5$ is fixed and $k$ is large, then $r(K_p,K_k)\ge \Omega(k^{\frac{p+1}{2}}(\log k)^{\frac{1}{p-2}-\frac{p+1}{2}})$ (which is also a lower bound for $r(K_p,B_{k,t})$) by Bohman and Keevash \cite{bk}, improving the best known lower bound due to Spencer \cite{spe} by a factor $(\log k)^{\frac{1}{p-2}}$.

\medskip
We also need the following stability result.
\begin{lemma}\label{stable}
  For any fixed graph $H$, integers $p\ge2$ and $b\ge1$, there exists $\delta>0$ such that the following holds for all
  $n\ge pb^2/\delta$. Let $a_1,\dots,a_p$ be positive integers  with $1=a_1\le a_2\leq \dots \leq a_{p-1}\leq b$ and $a_{p}\leq\delta n$, and let $G=K_p(1,a_2,\dots,a_{p})$.
For any graph $\Gamma$ on $N\ge (p-1)n|H|$ vertices containing no copy of $G$ and its complement $\overline{\Gamma}$ is $(K_1+nH)$-free,
there exists a partition $V(\Gamma)=\sqcup_{1\le i\le p-1} V_i$ such that
each vertex of $V_i$ has at most $a_2-1$ neighbors in $V_i$, for $i=1, 2, \dots, p-1$.

\medskip
Moreover, we may take $\delta$ with $0<\delta<\min\big\{\frac{1}{400a^{|H|+2}p^4},(100a^pp^{14p})^{-A}\big\}$, where $a=\sum_{i=1}^{p-1}a_i$ and $A=\prod_{i=1}^{p-1}a_i$.
\end{lemma}
{\bf Proof.} We assume $p\ge3$ since the assertion is trivial for $p=2$. Let $h$ denote the number of vertices in   $H$ and let
$a=\sum_{i=1}^{p-1}a_i$, and $A=\prod_{i=1}^{p-1}a_i.$  Let $\varepsilon=1/(100a^2p^4)$.
We may choose $\delta$ with $0<\delta<\min\big\{\frac{1}{400a^{h+2}p^4},(100a^pp^{14p})^{-A}\big\}$. We follow the notation in Lemma \ref{stability} to select $\gamma$ and $\eta$ such that $2(a+1)^h\delta\le\gamma\le\frac{\varepsilon}{2}$ and $\eta=\delta^{1/A}<p^{-10p}\varepsilon$. Furthermore, we assume that
\begin{align}
\label{lowerbound}
n\ge pb^2/\delta\ge\max\left\{b^2/\delta, a/\delta\right\},
\end{align}
and let $$\ell=r(K_h,G).$$

Note that $G$ is a subgraph of the book graph $B_{a,a+\lceil\delta n\rceil}$, it follows from Lemma  \ref{fox} that
\begin{align}\label{ell}
\ell =r(K_h, G) \leq r(K_h, B_{a,a+\lceil\delta n\rceil}) \le (a+1)^h (a+\lceil\delta n\rceil)<2(a+1)^h\delta n.
\end{align}
Since $\overline{\Gamma}$ contains no copy of $K_1+nK_h$, by a similar argument as Lemma \ref{deg}, we obtain that $d_{\overline{\Gamma}}(v)<r(G,nK_h)\le \ell+(n-1)h$ for any vertex $v$ in $\Gamma$.
 Therefore for any vertex $v$ in $\Gamma$,
\begin{align}\label{deg-low-2}
d_\Gamma(v)= N-1-d_{\overline{\Gamma}}(v)\ge N-\ell-(n-1)h\ge\left(1-\frac{1}{p-1}-\gamma\right)N
\end{align}
using $\gamma\ge2(a+1)^h\delta$.
Moreover, from Lemma \ref{Erdos},  $\Gamma$ has at most $\eta N^p$ copies of $K_p$   since $\Gamma$ is $G$-free. It follows from  Lemma \ref{stability} that there is a partition $V(\Gamma)=\sqcup_{1\le i\le p-1} V_i$ such that (i)-(iv) of Lemma \ref{stability} hold.

Let $z$ denote the total number of non-edges of $\Gamma$. From Lemma \ref{stability} (iii),
\[
z\le {p-1\choose 2}p^2\varepsilon|V_i||V_j|<p^2\varepsilon N^2.
\]
Let $s=\lceil p\sqrt{\varepsilon}N\rceil$. Clearly, $z< s^2$.
 Suppose that some vertex $v\in V_1$ satisfies $d_{V_1}(v)\ge2s$. Then, Lemma \ref{stability} (iv) implies that $d_{V_i}(v)\ge 2s$ for $2\le i\le p-1$. Let $U_i$ denote the neighborhood of $v$ in $V_i$. Then, by Lemma \ref{turan}, the subgraph of  $\Gamma$ induced by $\sqcup_{1\le i\le p-1} U_i$ contains at least $$(2s)^{p-3}(4s^2-z)\ge3\cdot 2^{p-3}s^{p-1}>(2p)^{p-2}(\sqrt{\varepsilon})^{p-1}N^{p-1}>\delta^{1/A}N^{p-1}\ge\delta^{1/(a_2\cdots a_{p-1})}N^{p-1}$$
distinct copies of $K_{p-1}$. Thus, by Lemma \ref{Erdos}, the neighborhood of $v$ contains a copy of $K_{p-1}(a_2,\dots,a_p)$. This leads to  a contradiction to the fact that $\Gamma$ contains no $K_p(1,a_2,\dots,a_p)$.
Therefore, $d_{V_1}(v)<2s$ for each $v\in V_1$.  Similarly, for any  $2\le i\le p-1$ and $u \in V_i$, we have
 $d_{V_i}(u)<2s$.

Suppose that some vertex $v\in V_i$ satisfies $d_{V_j}(v)\le(1-4p^2\sqrt{\varepsilon})|V_j|$ for some $j\neq i$. From the above, the vertex $v$ has at least $|V_i|-2p\sqrt{\varepsilon}N-1$ non-neighbors in $V_i$. Thus the total number of non-neighbors of $v$ is at least $4p^2\sqrt{\varepsilon}|V_j|+|V_i|-2p\sqrt{\varepsilon}N-1$, which is at least
\[
\left(1+4p^2\sqrt{\varepsilon}\right)\left(\frac{N}{p-1}-\sqrt{2\varepsilon}N\right)-2p\sqrt{\varepsilon}N-1>\left(\frac{1}{p-1}+\varepsilon\right)N,
\]
since $\varepsilon=1/(100a^2p^4)$.  By noting $\gamma\le \varepsilon/2$, we have  a contradiction to the fact that the minimum degree of $\Gamma$ is at least $(1-\frac{1}{p-1}-\gamma)N$ from (\ref{deg-low-2}). Therefore, for each vertex $v\in V_i$, $1\le i\le p-1$, we have
\begin{align}\label{ot-d}
\text{ $d_{V_j}(v)>(1-4p^2\sqrt{\varepsilon})|V_j|$  for any $j\neq i$.}
\end{align}

Now, suppose to the contrary that the assertion of the lemma does not hold. Without loss of generality, we may assume that there exists some vertex $v_1\in V_1$ such that $v_1$ has $a_2$ neighbors in $V_1$. It follows from  (\ref{ot-d}) that the vertex $v_1$ and these $a_2$ neighbors must have at least $(1-4p^2\sqrt{\varepsilon}(a_2+1))|V_2|
$ common neighbors in $V_2$.
Note that
\begin{align}\label{det}
|V_{i}|-a\cdot 4p^2\sqrt{\varepsilon}|V_{i}|
\ge(1-4a p^2\sqrt{\varepsilon})\left(\frac{N}{p-1}-\sqrt{2\varepsilon}N\right)
>hn/2>\delta n >a_3
\end{align}
by using the facts that $a=\sum_{i=1}^{p-1}a_i$, $\delta< 1/2$,  $\varepsilon=1/(100a^2p^4)$, and
 $\big||V_i|-\frac{N}{p-1}\big|\leq\sqrt{2\varepsilon}N$ for each $i\in[p-1]$, from Lemma \ref{stability} (ii).
 Therefore,
 the vertex $v_1$ and its $a_2$ neighbors in $V_i$ must have at least $a_3$ common neighbors in $V_3$. We can then inductively apply (\ref{ot-d}) and   (\ref{det}) to obtain a copy of $K_p(1,a_2,\dots,a_{p})$ in $\Gamma$, which leads to a contradiction. The assertion is proved. \hfill$\Box$

\medskip\noindent
{\bf Proof  of Theorem \ref{main-2}.}
The lower bound follows from Theorem \ref{low-bou}. It suffices to establish the upper bounds. We rely heavily on Lemma \ref{stable} and we follow all the definitions in its proof. In particular, we choose $n\ge pb^2/\delta$ where $0<\delta<\min\big\{\frac{1}{400a^{h+2}p^4},(100a^pp^{14p})^{-A}\big\}$. Recall  $G=K_p(a_1,\dots,a_{p})$.



\medskip\noindent
{\it Case 1:} \ Either $nh+a_2-1$  is even or $a_2-1$ is even.

\medskip
For this case, let $N_1=(p-1)(hn+a_2-1)+1$. Suppose on the contrary that there exists a graph $\Gamma$ on $N_1$ vertices such that $\Gamma$ is $G$-free and its complement $\overline{\Gamma}$ contains no copy of $K_1+nH$.
From Lemma \ref{stable}, there exists a partition $V(\Gamma)=\sqcup_{1\le i\le p-1} V_i$ such that the following holds:

\medskip
($*$) {\em Each vertex of $V_i$ has at most $a_2-1$ neighbors in $V_i$, for $i=1, 2, \dots, p-1$.}

\medskip
We may assume that $V_1$ is the largest part among $V_1,\dots,V_{p-1}$. Thus $|V_1|\ge\lceil \frac N{p-1}\rceil\ge hn+a_2.$
From ($*$), any vertex $x\in V_1$ has at most $a_2-1$ neighbors in $V_1$.
Therefore, there exists an independent set $W\subset V_1$ with
\[
|W|\ge |V_1|/a_2 \ge hn/a_2+1,
\]
and any vertex in $W$ has at least $(hn+a_2-1)-(a_2-1)=hn$ non-neighbors in $V_1$. Fix a vertex $w\in W$, and let $X$ be the non-neighborhood of $w$ in $V_1$. Clearly, $|X|\ge hn$.

\begin{claim}\label{nKh}
 $X$ contains $n$ disjoint independent sets of size $h$.
\end{claim}
{\bf Proof.}  Set $|X\setminus W|=\ell+n_1h+h'$, where $0\le
h'<h$ and $\ell=r(K_p(a_1,\dots,a_{p}),K_h)$. Since there is no $K_p(a_1,\dots,a_{p})$, we then can find at least
$n_1+1$ disjoint independent sets of size $h$ in $X\setminus W$. Let $X_0$ denote the remaining vertices in $X\setminus W$  by deleting the vertices of these disjoint independent sets of size $h$. Then $|X_0|<\ell<2(a+1)^h\delta n$ from (\ref{ell}). From ($*$), each vertex $x\in X_0$ has at most $a_2-1$ neighbors in $W\setminus\{w\}$ and therefore,  there is a subset $W_1$ of $W\setminus \{w\}$
consisting of vertices non-adjacent to any vertex in $X_0$ satisfying
$$
|W_1|\ge |W\setminus \{w\}|-(a_2-1)\ell\ge hn/a_2-2a_2(a+1)^h\delta n\ge(h-1)\ell
$$
where the last inequality follows from the fact that $\delta$ is sufficiently small. Thus any vertex $x\in X_0$ and $h-1$ vertices of $W_1$ form an
independent set of size $h$. Let $W_0\subset W_1$ be the set  consisting of the vertices that have been accounted for.  The remaining vertices in $W\setminus W_0$ clearly forms an independent set. Since $|X|\ge hn$, we can definitely obtain $n$ disjoint independent sets of size $h$  as desired. \hfill$\Box$

\medskip
Claim \ref{nKh} implies that $X\cup\{w\}$ yields a copy of $K_1+nK_h$ in the complement of $\Gamma$ with center $w$. This leads to a contradiction.

\medskip\noindent
{\it Case 2:} \ Both $nh+a_2-1$  and $a_2-1$ are odd.

\medskip
For this case, let $N_2=(p-1)(hn+a_2-2)+1$. Suppose on the contrary that there exists a graph $\Gamma$ on $N_2$ vertices such that $\Gamma$ is $G$-free and its complement $\overline{\Gamma}$ contains no copy of $K_1+nH$. We shall show that this will lead to a contradiction.
From Lemma \ref{stable}, there exists a partition $V(\Gamma)=\sqcup_{1\le i\le p-1} V_i$ such that
each vertex of $V_i$ has at most $a_2-1$ neighbors in $V_i$, for $i=1, 2, \dots, p-1$.

If there exists some part $V_i$ has size at least $hn+a_2$, then we are done by a similar argument as in Case 1. So we may assume that $|V_i|\le hn+a_2-1$ for $1\le i\le p-1$. There must exist one part, say $V_1$, of size  at least $\lceil\frac{N_2}{p-1}\rceil=hn+a_2-1$. Thus,  $|V_1|=hn+a_2-1$.
\begin{claim}\label{a2-1-rg}
Each vertex of $V_1$ has exactly $a_2-1$ neighbors in $V_1$.
\end{claim}
{\bf Proof.}  Suppose, on the contrary, there is a vertex $w\in V_1$ that has at most $a_2-2$ neighbors in $V_1$. Then it has at least $hn$ non-neighbors in $V_1$. Let $X$ be the non-neighborhood of $w$ in $V_1$. Clearly, $|X|\ge hn$.
Since any vertex $x\in X$ has at most $a_2-1$ neighbors in $X$,  there is an independent set $W\subset X$ with $|W|\ge |X|/a_2 \ge hn/a_2.$ By a similar argument as in Claim \ref{nKh}, $X$ contains $n$ disjoint independent sets of size $h$, which together with $w$ yield a copy of $K_1+nK_h$ in the complement of $\Gamma$ which is impossible.\hfill$\Box$

\medskip
From Claim \ref{a2-1-rg},  we conclude that the subgraph of $\Gamma$ induced by $V_1$ is $(a_2-1)$-regular. However, such a subgraph of order $hn+a_2-1$ does not exist since both $hn+a_2-1$ and $a_2-1$ are odd.
This  completes the proof of  Theorem \ref{main-2}.
\hfill$\Box$

\section{Problems and remarks}\label{clu}

As a central subject in combinatorics, the problem of determining the exact values of Ramsey numbers is notoriously difficult. The study of goodness of Ramsey numbers follows an opposite path, in search of graphs that can achieve the  (relatively weak) lower bounds or with small discrepancies. The main results in this paper is along this line of approaches. Nevertheless, numerous questions remain, some of which we mention here.

\begin{problem}
Find a characterization for graphs $H$ that is $K_p$-good. Namely, determine the family of graphs $H$ satisfying $r(K_p, H)=(p-1)(|H|-1) +1$.

{\rm So far, it is known that this family includes connected graphs with bounded maximum degree and  small bandwidth \cite{abs}, connected graphs with bounded degeneracy satisfying certain locally sparse conditions \cite{nr09}, etc., but the list is  far from complete.
We remark that trees are included in the above list (as seen in  (\ref{chvatal})), belonging to the family of bounded degeneracy.
The {\it degeneracy} $d(H)$ of a graph $H$ is the smallest natural
number $d$ such that every induced subgraph of $H$ has a vertex of degree at most $d$. For example, a tree has degeneracy $1$.
}
\end{problem}
\begin{problem}
Let $p\ge3$ and $a_3,\dots,a_p$ be integers  with $2\le a_3\leq \dots \leq a_{p}$. Find a characterization for graphs $H$ that is $K_p(1,2,a_3, \dots, a_p)$-good.

 {\rm The main theorems in this paper provide some  hints  in this direction. }
\end{problem}

\begin{problem}
Give some classifications for graphs $H$ with low discrepancies from $G$-goodness.

 {\rm Of course, this problem may be too general or too hard to tackle. Here we just intend to point out numerous possible directions.}
 
 \end{problem}

\begin{problem}
{\rm Corollary \ref{fan} shows that if $n\ge Cp^2$,  then $F_n$ is $K_p$-good where $C\approx 26.228$.}

 It would be interesting to improve the lower bound of $n$ further, e.g., is it true for $n\ge \Omega(p)$?
 Moreover, it would be interesting to improve the lower bounds of $n$ in Theorem \ref{main} and Theorem \ref{main-2}.

\end{problem}
\begin{problem}
{\rm Nikiforov and Rousseau \cite[Theorem 2.1]{nr09} established an extremely general Ramsey goodness theorem for several families of graphs. However, the quantitative dependence between the graph sizes involved
are tower-type since the proofs rely on Szemer\'edi's regularity lemma \cite{sz}. Fox, He and Wigderson \cite{fox} showed that for every $k, p, t\geq 2$, there exists $\delta>0$ such that the following holds for all large $n$. Let $1\leq a_1\leq \dots \leq a_{p-1}\leq t$ and $a_{p}\leq\delta n$ be positive integers. If $a_1=a_2=1$, then $r(K_p(a_1,\dots,a_{p}), B_{k,n})= (p-1)(n-1)+1$. Their proof does not use the regularity lemma, and thus double-exponential bounds on $\delta$ will suffice. In particular, $B_{k,n}$ is $K_p$-good provided $n\ge2^{k^{10p}}$. We refer the reader to \cite{fl} for some extended results.}

Fox, He and Wigderson \cite{fox} asked if it is possible to completely eliminate the use of the regularity lemma from the proof of \cite[Theorem 2.1]{nr09}, which would likely lead to superior quantitative bounds.

\end{problem}

\bigskip
{\bf Acknowledgment.}
We are grateful to the anonymous referee for giving invaluable comments and suggestions which  greatly improve the presentation of the manuscript (especially Theorem \ref{main} and Theorem \ref{low-bou}).

\end{spacing}

\end{document}